\documentclass[12pt]{amsart}
\usepackage[margin=1in]{geometry}
\usepackage{times}
\usepackage{amsthm,amssymb}
\usepackage{amsmath,xypic}
\usepackage[usenames,dvipsnames]{color}
\usepackage[
colorlinks=true,hyperindex, citecolor=cyan, urlcolor=cyan]{hyperref}
\usepackage{epigraph}
\usepackage{graphicx}
\usepackage{fancyhdr} 
\usepackage{enumitem}
\newcommand{\be}{\begin{equation}}
\newcommand{\ee}{\end{equation}}
\newcommand{\bes}{\begin{equation*}}
\newcommand{\ees}{\end{equation*}}
\newcommand{\bea}{\begin{eqnarray}}
\newcommand{\eea}{\end{eqnarray}}
\newcommand{\beas}{\begin{eqnarray}}
\newcommand{\eeas}{\end{eqnarray}}
\newcommand{\ben}{\begin{note}}
\newcommand{\een}{\end{note}}
\newcommand{\bexl}{\vskip0.1em\noindent\hrulefill\vskip1em\begin{ExerciseList}}
\newcommand{\eexl}{\end{ExerciseList}\hrulefill}

\newcommand{\bthm}{\begin{theorem}}
\newcommand{\ethm}{\end{theorem}}
\newcommand{\bpro}{\begin{prop}}
\newcommand{\epro}{\end{prop}}
\newcommand{\bcor}{\begin{corollary}}
\newcommand{\ecor}{\end{corollary}}
\newcommand{\bcon}{\begin{conjecture}}
\newcommand{\econ}{\end{conjecture}}
\newcommand{\bp}{\begin{proof}}
\newcommand{\ep}{\end{proof}}
\newcommand{\blem}{\begin{lemma}}
\newcommand{\elem}{\end{lemma}}
\newcommand{\bn}{\begin{note}}
\newcommand{\en}{\end{note}}
\newcommand{\benum}{\begin{enumerate}}
\newcommand{\eenum}{\end{enumerate}}
\newcommand{\bed}{\begin{defn}}
\newcommand{\eed}{\end{defn}}
\newcommand{\brem}{\begin{remark}}
\newcommand{\erem}{\end{remark}}

\newcommand{\btik}{\begin{tikzpicture}\begin{axis}[scale=0.5,axis y line=center, axis x line=middle]}
\newcommand{\etik}{\end{axis}\end{tikzpicture}}

\let\into=\hookrightarrow

\newcommand{\upperRomannumeral}[1]{\uppercase\expandafter{\romannumeral#1}}

\newtheorem{theorem}[equation]{Theorem}      
\newtheorem{lemma}[equation]{Lemma}          %
\newtheorem{corollary}[equation]{Corollary}  
\newtheorem{proposition}[equation]{Proposition}

\theoremstyle{definition}

\theoremstyle{definition}
\newtheorem{defn}[equation]{Definition}
\theoremstyle{remark}

\theoremstyle{definition}
\newtheorem{remark}[equation]{Remark}

\numberwithin{equation}{section}

\let\congruent=\equiv
\let\into=\hookrightarrow
\let\isom=\simeq

\let\tensor=\otimes

\newcommand{\Pic}{{\rm Pic\,}}

\newcommand{\Z}{{\mathbb Z}}

\renewcommand{\int}{\operatorname{int}}
\renewcommand{\O}{{\mathcal O}}
\renewcommand{\P}{{\mathbb P}}

\renewcommand{\bpro}{\begin{proposition}}
\renewcommand{\epro}{\end{proposition}}

\let\citep=\cite 

\begin{document}

\title[]{A remark on Ulrich and ACM bundles}%
\author{Kirti Joshi}%
\address{Math. department, University of Arizona, 617 N Santa Rita, Tucson
85721-0089, USA.} \email{kirti@math.arizona.edu}

\thanks{}%


\begin{abstract}
I show that on any smooth, projective ordinary curve of genus at least two and a projective embedding,  there is a natural example of a stable Ulrich bundle for this embedding: namely the sheaf $B^1_X$ of locally exact differentials twisted by $\O_X(1)$ given by this embedding and in particular there exist ordinary varieties of any dimension which carry Ulrich bundles. In higher dimensions, assuming $X$ is Frobenius split variety I show  that $B^1_X$ is an ACM bundle and if $X$ is also a Calabi-Yau variety and $p>2$ then $B^1_X$ is not a direct sum of line bundles. In particular I show that $B^1_X$ is an ACM bundle on any ordinary Calabi-Yau variety. I also prove a characterization  of  projective varieties with trivial canonical bundle such that $B^1_X$ is ACM (for some projective embedding datum): all such varieties  are  Frobenius split  (with trivial canonical bundle).
\end{abstract}
\maketitle
\epigraph{My eyes are used to sights like these:\\ I stand between familiar trees.}{Theodore Roethke (\cite{roethke38})}

\tableofcontents
\section{Introduction} Let $X$ be a smooth projective variety over an algebraically closed field $k$ equipped with a projective embedding $X\into \P^n$ and $\O_X(1)$ the very ample line bundle on $X$ provided by this embedding.  Let $E$ be a vector bundle on $X$. Then $E$ is an \emph{Ulrich bundle} if $H^i(X,E(-i))=0$ for $i>0$ and $H^j(X,E(-j-1))=0$ for $j<\dim(X)$ (see \cite[Proposition 2.1(b)]{eisenbud03}). It has been conjectured that Ulrich bundles exists on every projective variety. This is known in very small number of cases (see (\cite{beauville17}) for an excellent survey of this topic). In general the problem of constructing Ulrich bundles is difficult and many constructions are specific to the sort of variety under consideration and often provide Ulrich bundles small rank. 

Now suppose $k$ has characteristic $p>0$. The purpose of this note is to record the following elementary remark which provides a natural example of a stable Ulrich bundle of large rank on ordinary curves. A smooth projective curve  $X$ over an algebraically closed field is \emph{ordinary} if and only if Frobenius map $H^i(X,\O_X)\to H^i(X,\O_X)$ is an isomorphism for all $i\geq 0$. Let $d:\O_X\to \Omega^1_X$ be the differential. The image $B_X^1=d(\O_X)$ is a subsheaf of $\Omega^1_X$ consisting of locally exact differentials. As $d(f^pg)=f^pdg$ one sees that $B^1_X$ is a locally free subsheaf $B^1_X\subset F_*(\Omega^1_X)$ and  one has the fundamental exact sequence
\be\label{eq:fundamental-seq} 
0\to \O_X\to F_*(\O_X)\to B^1_X\to 0. 
\ee
Thus $B^1_X$ is locally free of rank $p-1$ and degree $(p-1)(g-1)$. 

My remark is this: if $X$ is a smooth, projective curve of any genus then $E=B^1_X(1)$ is an Ulrich bundle of rank $p-1$ for any smooth, projective embedding $X\into \P^n$ and $\O_X(1)$ the very ample line bundle provided by this embedding if and only if $X$ is ordinary (note that by \cite{joshi04}, $B^1_X$ is a stable bundle for $g\geq 2$). This also provides examples of ordinary varieties of any dimension equipped with Ulrich bundles (see Theorem~\ref{th:main}) and in Theorem~\ref{th:more-examples} I provide many more examples of Ulrich bundles on curves in characteristic $p$. The bundles $B^1_X$ plays an important role in the algebraic geometry of smooth, projective curves: for instance it plays a central role in the work of Michel Raynaud (see \cite{raynaud00}) and Akio Tamagawa on the fundamental groups of smooth projective curves, notably in Tamagawa's proof of Raynaud's conjecture (see \cite{tamagawa04}). My remark shows that if the curve is ordinary the syzygies of this bundle are linear  (in the sense of \cite[Proposition 2.1(d)]{eisenbud03}) and suggests that in characteristic $p$ the arithmetic nature of the curve is  mirrored in the structure of syzygy tables of bundles on the curve--a connection which is perhaps not manifestly obvious. In general i.e. if $X$ is not ordinary, the result is less explicit: there always exists a line bundle $L$ of degree zero such that $B^1_X\tensor L(1)$ is an Ulrich bundle on $X$. 

Recall that a vector bundle $E$ on a smooth, projective variety $X$ is an \emph{ACM bundle (arithmetically Cohen-Macaulay bundle)} if $H^i(X,E(n))=0$ for all $0<i<\dim(X)$ and for all $n\in\Z$. In dimensions greater than one, any Ulrich bundle is an ACM bundle. I do not know if a suitable twist  of $B^1_X$ is an Ulrich bundle in higher dimensions. In Theorem~\ref{th:fsplit} I show that $B^1_X$ is an ACM bundle on any Frobenius split Calabi-Yau variety (for any embedding of $X$) and in Corollary~\ref{co:calabi-yau} I show that if $p>2$, $X$ is Calabi-Yau and Frobenius split then $B^1_X$ is an ACM bundle which is not a direct sum of line bundles. Finally let me remark that if $X$ is a smooth, projective Calabi-Yau variety in characteristic $p>0$ and $X$ satisfies Kodaira vanishing then $F_*(\O_X)$ is always an ACM bundle on $X$ for any projective embedding of $X$. 

It is a pleasure to thank N.~Mohan Kumar for correspondence and a number of comments. I would also like to thank the referee for a number of suggestions in particular I have included Section~\ref{sec:examples} which provides examples of varieties satisfying various hypothesis of my theorems at the referee's request.

\section{Ulrich bundles on curves}
\newcommand{\sE}{E}
\newcommand{\sL}{L}
\newcommand{\boxten}{\mathop\boxtimes\displaylimits}
For the rest of the paper let $k$ be an algebraically closed field of characteristic $p>0$. For $k$-schemes $X_1$ (resp. $X_2$) and sheaves $E$ (resp. $F$) on $X_1$ (resp. $X_2$) let $E\boxten F$ be the tensor product $q_1^*(E)\tensor q_2^*(F)$ where $q_1$ (resp. $q_2$) is the projection from $X_1\times X_2$ to the factor $X_1$ (resp. factor $X_2$). 
\bthm\label{th:main}
Let $k$ be an algebraically closed field of characteristic $p>0$. Let $(X_i,\O_{X_i}(1))$ for $i=1,\ldots,m$ be smooth, projective ordinary curves of genus at least two and each equipped with projective embedding $q_i:X_i\into\P^{n_i}$ and $\O_{X_i}(1)$ the very ample line bundle given by this embedding. Then the bundle 
$$\sE = \boxten_{i=1}^m  B^1_{X_i}(i)$$ 
is an Ulrich bundle on $Y=X_1\times X_2\times \cdots \times X_m$ with respect to the Segre embedding given by $\sL_m=\boxten_{i=1}^m \O_{X_i}(1)$
\ethm
Let me remark that for all $g\geq 1$ there is a dense open subset of the moduli of smooth, proper curves of genus $g$ which parameterizes ordinary curves. So ordinarity is a genericity condition on curves of genus $g\geq 1$.
\bp 
Let me first prove this for $m=1$. So I have to prove that if $X$ is a smooth, projective and ordinary curve embedded in projective space $X\into \P^n$ and $\O_X(1)$ is the very ample line bundle given by this embedding then $\sE = B^1_X(1)$ an Ulrich bundle. To prove this it suffices to observe that $H^j(X,E(-1))=H^j(X,B^1_X)=0$ for all $j\geq 0$ as $X$ is ordinary. Thus $B^1_X(1)$ is an Ulrich bundle on $X$ for any smooth, projective embedding $X\into \P^n$. Now I prove the general claim by induction on $m$. Clearly the result is true for $m=1$. Suppose the assertion  is true for some $m-1\geq 1$. Then $\boxten_{i=1}^{m-1}  B^1_{X_i}(i)$ is an Ulrich bundle on $X_1\times \cdots \times X_{m-1}$ for $\sL_{m-1}$. Now by (\cite[Proposition 2.6]{eisenbud03}) 
$$\boxten_{i=1}^{m-1}  B^1_{X_i}(i)\boxten B^1_{X_m}(1)(m-1)=\boxten_{i=1}^m  B^1_{X_i}(i) $$ is Ulrich bundle on $X_1\times \cdots \times X_m$ for $\sL_m=\sL_{m-1}\boxten\O_{X_m}(1)$. This proves the theorem.
\ep

\bcor 
Let $X\into \P^n$ be a smooth, projective  curve over an algebraically closed field $k$ of characteristic $p>0$ and $\O_X(1)$ the very ample line bundle provided by this embedding. 
\benum[label={\bf{(\arabic{*})}}]
\item Then $B^1_X(1)$ is an Ulrich bundle (of rank $p-1$ and degree $(p-1)(g-1)$) if and only if $X$ is ordinary.
\item In particular if $X$ is ordinary then $X$ carries a canonical, stable Ulrich bundle: $B^1_X(1)$.
\item There always exists a line bundle $L$ of degree zero on $X$ such that $B^1_X\tensor L(1)$ is a stable Ulrich bundle on $X$. 
\eenum
\ecor

\bp
The rank and degree calculations are well-known \cite{raynaud82} and stability of $B^1_X$ (ordinarity is not needed for this) is due to (\cite{joshi04}). The ordinary case is immediate from definitions (of Ulrich bundle and of ordinarity), indeed observe that  $H^i(X,B^1_X(1)(-1))=H^i(X,B^1_X)$ for all $i\geq 0$ if and only if $X$ is ordinary, the first assertion is immediate. One sees from \cite{raynaud82} that there exists a dense open subset in $\Pic(X)$ consisting of line bundles $L$ of degree zero such that $H^*(X,B^1_X\tensor L)=0$. Then $E= B^1_X\tensor L(1)$ is a stable Ulrich bundle on $X$.
\ep

The following result provides many more examples of Ulrich bundles on any smooth, projective curve.
\bthm\label{th:more-examples} 
Let $X\into\P^n$ be a smooth, projective curve equipped with a projective embedding and $\O_X(1)$ the very ample line bundle given by this embedding. Suppose $X$ is ordinary. Then there is a non-empty dense open subset $\mathcal{SU}(r,\O_X)$ of moduli space of semistable bundles $V$ of rank $r$ and trivial determinant such that $E=B^1_X\tensor V(1)$ is an Ulrich bundle.
\ethm

\bp 
It is sufficient to produce a dense open set of semistable bundles $V$ of rank $r$ such that $B^1_X\tensor V$ has no cohomology. Since this cohomology vanishing is an open condition one may simply note that the trivial bundle of rank $r$ has this property $H^i(X,B^1_X\tensor \O_X^{\oplus r})=0$ for $i\geq 0$. This is true by ordinarity of $X$. This proves the assertion.
\ep

\section{ACM bundles on Frobenius split varieties}
I do not know if $B^1_X$ is an Ulrich bundle in higher dimensions and this is unlikely to hold in complete generality even for ordinary varieties. If $\dim(X)\geq 2$ then one sees from (\cite{coskun13}) that $B^1_X$ is an Ulrich bundle if and only if $B^1_X$ is ACM bundle and its Hilbert polynomial satisfies some condition. The Hilbert polynomial condition does sometimes fail to hold (see Remark~\ref{re:quartic}) below. In what follows I will write $\omega_X$ for the canonical line bundle of $X$.

Following (\cite{mehta85}) I say that a smooth, projective variety $X$ is \emph{Frobenius split} if \eqref{eq:fundamental-seq} splits as a sequence of $\O_X$-modules. For a Frobenius split variety one has $H^i(X,B_X^1)=0$ for all $i\geq 0$ (see \cite{joshi03}). 

Recall that a smooth, projective variety $X$ is \emph{ordinary} in the sense of (\cite{illusie83b}) if $H^j(X,d(\Omega_X^{i}))=0$ for all $i,j\geq 0$ (here $d:\Omega^i_X\to\Omega_X^{i+1}$ is the usual differential). Note that any smooth, proper curve is ordinary in this sense if and only if it is ordinary in the earlier sense).

For a smooth, projective variety $X$ and an ample line bundle $L$ on $X$, I say that $X$ satisfies \emph{Kodaira vanishing  for $L$} if $H^i(X,L^n)=0$ for all $i<\dim(X)$ and all integers $n<0$ and that $X$ satisfies Kodaira vanishing if $X$ satisfies \emph{Kodaira vanishing} for any ample line bundle on $X$. 

It was shown in (\cite{mehta85}) that Frobenius split varieties satisfy Kodaira vanishing theorem. Moreover it was also shown in  loc. cit.  that Frobenius split varieties  also satisfy the following vanishing: $H^i(X,L)=0$ for any ample line bundle and any integer $i>0$. This will be used in proof of Theorem~\ref{th:fsplit0}.

Theorem~\ref{th:fsplit0} shows that for any Frobenius split smooth, projective variety $B^1_X$ is an ACM bundle. Theorem~\ref{th:fsplit} provides a characterization of  varieties with trivial canonical bundle such that  $B^1_X$ is  an ACM bundle (for some embedding datum $(X\into\P^n,\O_X(1))$): all such varieties are Frobenius split.  By (\cite[Proposition~3.1]{joshi03}) Frobenius split varieties with trivial canonical bundles include ordinary varieties with trivial canonical bundle and hence Theorem~\ref{th:fsplit}  includes as special cases: ordinary abelian varieties, ordinary Calabi-Yau varieties. In Theorem~\ref{co:calabi-yau} I show that if $X$ a Calabi-Yau variety and $p>2$ then $B^1_X$ is an ACM bundle which is not a direct sum of line bundles. 

For examples of varieties satisfying hypothesis of this and subsequent theorems see Section~\ref{sec:examples}.

\bthm\label{th:fsplit0}
Let $X/k$ be a smooth, projective and Frobenius split variety and  equipped with  $\O_X(1)$ given by a projective embedding $X\into\P^n$. Then $B^1_X$ is an ACM bundle on $X$. 
\ethm
\bp 
By \eqref{eq:fundamental-seq} and (\cite[Proof of Theorem~4.2]{joshi03}) one has $H^i(X,B^1_X)=0$ for all $i\geq 0$. Hence it remains to prove that for any integer $m\neq0$ one has $H^i(X,B^1_X(m))=0$ for $1\leq i<\dim(X)$. By the well-known result of (\cite{mehta85}) one sees that the $X$ satisfies Kodaira vanishing for all ample line bundles and also $H^i(X,L)=0$ for all $i>0$ for any ample line bundle $L$ on $X$. Hence one has in particular that $H^i(X,\O_X(m))=0$ for all $m>0$ and $1\leq i<\dim(X)$. As $X$ satisfies Kodaira vanishing one has $H^i(X,\O_X(m))=0$ for $1\leq i<\dim(X)$ for all integers $m<0$. As $B_X^1$ is a direct summand of $F_*(\O_X)$ (by Frobenius splitting hypothesis) one has $B^1_X(m)$ is a direct summand of $F_*(\O_X)(m)$ and hence $H^i(X,B^1_X(m))=0$ for all $m\neq0$ and $1\leq i<\dim(X)$ as this cohomology is a direct summand of $H^i(X,F_*(\O_X)(m))=H^i(X,\O_X(pm))=0$ for $1\leq i<\dim(X)$.
\ep

The following special case is useful:
\bcor
Let $X\into\P^n$ be a smooth, projective, ordinary abelian variety or an ordinary Calabi-Yau variety equipped with $\O_X(1)$ given by this embedding.
Then $B^1_X$ is an ACM bundle on $X$.
\ecor
\bp
To prove this it suffices to note that by (\cite[Proposition~3.1(b)]{joshi03}) any abelian or Calabi-Yau variety which is also ordinary in the sense of (\cite{illusie83b}) is also Frobenius split. So the result is immediate from Theorem~\ref{th:fsplit0}. 
\ep

The next theorem characterizes varieties with trivial canonical bundles such that $B^1_X$ is an ACM bundle: these are precisely Frobenius split varieties with trivial canonical bundle (by (\cite{joshi03})  ordinary varieties with trivial canonical bundle are also Frobenius split).

\bthm\label{th:fsplit}
Let $X/k$ be a smooth, projective  variety satisfying $\omega_X=\O_X$. 
Then the following conditions are equivalent.
\benum[label={{\bf(\arabic{*})}}]
\item\label{th-fsplit1} $X$ is Frobenius split.
\item\label{th-fsplit2} For any embedding datum $(X\into\P^n,\O_X(1))$, $B^1_X$ is an ACM bundle on $X$.
\item\label{th-fsplit3} For some embedding datum $(X\into\P^n,\O_X(1))$, $B^1_X$ is an ACM bundle on $X$.
\item\label{th-fsplit4} $H^{d-1}(X,B^1_X)=0$ with $d=\dim(X)$.
\eenum
\ethm

Before giving the proof of Theorem~\ref{th:fsplit}, let me record the following simple observation which may be of independent interest.
\bpro\label{pr:push-acm} 
Let $X\into\P^n$ be a smooth, projective variety equipped with $\O_X(1)$ given by this embedding. If $E$ is an ACM bundle on $X$, then so is $F_*(E)$. 
\epro

\bp 
Suppose $E$ is an ACM bundle for $X\into\P^n$. Then $H^i(X,E(m))=0$ for $0<i<\dim(X)$ and all $m\in\Z$. In particular $H^i(X,E(pm))=0$ for $0<i<\dim(X)$ and all $m\in \Z$. Then using the fact that Frobenius is a finite morphism and the projection formula $F_*(E(pm))=F_*(E)(m)$ one sees that $H^i(X,E(pm))=H^i(X,F_*(E)(m))=0$ for $0<i<\dim(X)$ and all $m\in \Z$. This proves the assertion.
\ep

\bthm\label{th:push-acm-2}
Suppose $X\into \P^n$ is a smooth, projective Calabi-Yau variety equipped with $\O_X(1)$ provided by this embedding. Suppose $X$ satisfies Kodaira vanishing for $\O_X(1)$. Then $F_*(\O_X)$ is an ACM bundle on $X$.
\ethm

\bp 
This follows from the previous proposition. As $X$ is Calabi-Yau $H^i(X,\O_X)=0$ for $1\leq i<\dim(X)$ and as $X$ satisfies Kodaira-vanishing for $\O_X(1)$ one sees that $H^i(X,\O_X(n))=0$ for all integers $n>0$ and $i>0$ and $H^j(X,\O_X(m))=0$ for $j<\dim(X)$ and all integers $m<0$; hence $H^i(X,\O_X(n))=0$ for all $n\in\Z$ and $1\leq i<\dim(X)$. In other words $\O_X$ is an ACM bundle on $X$ and so by Proposition~\ref{pr:push-acm} one sees that $F_*(\O_X)$ is also an ACM bundle on $X$. 
\ep

\bp[Proof of Theorem~\ref{th:fsplit}] 
The implication \ref{th-fsplit1}$\implies$\ref{th-fsplit2} is immediate from Theorem~\ref{th:fsplit0}.

The implication \ref{th-fsplit2}$\implies$\ref{th-fsplit3} is of course trivial.

The implication \ref{th-fsplit3}$\implies$\ref{th-fsplit4} is also clear: if $B^1_X$ is an ACM bundle for some embedding datum $(X\into\P^n,\O_X(1))$ then one has $H^i(X,B^1_X(m))=0$ for all integers $m$ and $1\leq i< d=\dim(X)$. Hence one has $H^{d-1}(X,B^1_X)=0$ which proves \ref{th-fsplit3}$\implies$\ref{th-fsplit4}.

So it remains to prove that \ref{th-fsplit4}$\implies$\ref{th-fsplit1}. So suppose $B^1_X$ is an ACM bundle. From the well-known criterion for Frobenius splitting (see \cite[Proposition 9]{mehta85b}) one knows that $X$ is Frobenius split if and only if the morphism $H^d(\omega_X)\to H^d(\omega_X^p)$ is injective. Since $\omega_X=\O_X$ and $B^1_X$ is an ACM bundle by hypothesis this follows from \eqref{eq:fundamental-seq} on taking cohomology:
$$
\xymatrix{
 \ar[r] & H^{d-1}(X,B^1_X)\ar@{=}[d]\ar[r]& H^d(X,\O_X)\ar[r]^F & H^d(X,\O_X)\ar [r]& H^d(X,B^1_X)\ar[r] &0.\\
&0&&&
}$$
This proves \ref{th-fsplit4}$\implies$\ref{th-fsplit1} and hence this completes the proof.
\ep

 Theorem~\ref{th:fsplit} has the following useful corollary.
\bcor\label{cor:fsplit-calabi-yau}
Let $X/k$ be a smooth, projective abelian or  Calabi-Yau  variety   equipped with  $\O_X(1)$ given by a projective embedding $X\into\P^n$. Then the following conditions are equivalent.
\benum[label={{\bf(\arabic{*})}}]
\item\label{th-fsplit-calabi-yau-1} $X$ is Frobenius split.
\item\label{th-fsplit-calabi-yau-2} $B^1_X$ is an ACM bundle on $X$ (for the given embedding $X\into \P^n$).
\item\label{th-fsplit-calabi-yau-3} $H^{d-1}(X,B^1_X)=0$ with $d=\dim(X)$.
\eenum
\ecor

For an ordinary abelian variety $B^1_X$ is a direct sum of line bundles (see \cite{sannai16}) on the other hand my next result shows that if $X$ is Calabi-Yau and $p>2$ then $B^1_X$ is an ACM bundle which is not a direct sum of line bundles.

\bthm\label{co:calabi-yau} 
Let $X\into \P^n$ be smooth, projective abelian or a Calabi-Yau variety with $\dim(X)\geq 2$ and equipped with $\O_X(1)$ given by this embedding. 
\benum[label={{\bf(\arabic{*})}}] 
\item If $X$ is a Frobenius split abelian variety (equivalently an ordinary abelian variety) then $B^1_X$ is a direct sum of ACM line bundles. In particular any ordinary abelian variety carries a distinguished finite set of ACM line bundles (for any given embedding datum).
\item If $X$ is a Frobenius split Calabi-Yau variety and $p>2$ then $B^1_X$ is an ACM bundle which not a direct sum of line bundles.
\eenum
\ethm

\bp 
By (\cite{sannai16}) one sees, as $X$ is a Frobenius split abelian variety (equivalently an ordinary abelian variety), that $B^1_X$ is a direct of sum line bundles and as $B^1_X$ is an ACM bundle, every line bundle direct summand of $B^1_X$ is also ACM and hence one has a distinguished finite set of ACM line bundles  on $X$ which is ACM (for any given embedding datum). This proves {(\bf1)}.

Let me prove {\bf(2)}. The only point which needs to be proved is that $B^1_X$ is not a direct sum of line bundles. Suppose $B^1_X$ is a direct sum of line bundles. Then so is $F_*(\O_X)=\O_X\oplus B^1_X$.  On the other hand if $p>2$, $F_*(\O_X)$ is a direct sum of line bundles, then as $\omega_X=\O_X$ is pseudo-effective, one sees by  (\cite[Theorem 5.5]{sannai17}) that $X$ is an abelian variety. But this is a contradiction as $X$ is Calabi-Yau variety. 
\ep
The preceding methods also show that
\bthm\label{th:fano} 
Suppose $X\into\P^n$ is a smooth, projective Fano variety which satisfies Kodaira vanishing. Then $F_*(\O_X)$ is an ACM bundle on $X$.
\ethm
\bp 
Recall that any Fano variety which satisfies Kodaira vanishing also satisfes $H^i(X,\O_X)=0$ for $1\leq i\leq \dim(X)$. By Proposition~\ref{pr:push-acm} it is enough to show that $\O_X$ is an ACM bundle. Since $X$ is Fano and satisfies Kodaira vanishing one has $H^i(X,\O_X(m))=0$ for all integers $m\leq 0$ and $1\leq i<\dim(X)$. Further for any $m>0$, by Serre duality and Kodaira vanishing $H^i(\O_X(m))\isom H^{n-i}(\omega_X\tensor\O_X(-m))=0$ for $1\leq i<\dim(X)$ as $\omega_X$ and $\O_X(-m)$ are both anti-ample line bundles. Hence $H^i(X,\O_X(m))=0$ for all $m\in\Z$ and $1\leq i<\dim(X)$ so $\O_X$ is an ACM bundle. Hence so is $F_*(\O_X)$.
\ep

\section{Examples and remarks}\label{sec:examples}
Let me provide example of varieties satisfying all the hypothesis of all these theorems. First of all Frobenius splitting, ordinarity are open conditions (on the base) in any flat family of Calabi-Yau varieties (see \cite{joshi03}). It is a well-known  (see (\cite{illusie90a})) that a general hypersurface of degree $\geq 1$ (and even complete intersections) in projective space are	ordinary. 

Now here is an explicit example of an ordinary Calabi-Yau variety. The Fermat Calabi-Yau hypersurface in $\P^n$:
$$x_0^{n+1}+x^{n+1}_1+\cdots +x_n^{n+1}=0$$
is ordinary if $p\congruent 1\bmod{n+1}$. This, together with (\cite{illusie90a}) provides examples of ordinary Calabi-Yau varieties satisfying hypothesis of Theorem~\ref{th:push-acm-2}.

Now let me provide an example of Frobenius split hypersurfaces in $\P^n$. By \cite{fedder87} or \cite[Proposition~6]{mehta85} if  $\deg(f)\leq n+1$ then $f(x_0,\ldots,x_n)=0$ is Frobenius split if and only if $f^{p-1}$ contains the term $a(x_0\cdots x_n)^{p-1}$ with $a\neq 0$. 

From this it is immediate that any hypersurface of  $\deg(f)<n+1$ in $\P^n$ is always Frobenius split. This provides examples of Fano varieties satisfying hypothesis of Theorem~\ref{th:fano}.

Now consider the most mathematically visible family of Calabi-Yau hypersurfaces in $\P^n$:
$$x_0^{n+1}+x^{n+1}_1+\cdots +x_n^{n+1}=\lambda x_0\cdot x_1\cdots x_n.$$
The above criterion shows that this is Frobenius split for general $\lambda$ if $p\congruent 1\bmod{n+1}$. This provides examples of Frobenius split, Calabi-Yau varieties satisfying hypothesis of Theorem~\ref{th:fsplit}.

Finally let me point out (\cite{joshi03}) as a reference for comparison of the two conditions: Frobenius splitting and ordinarity.

\brem 
Let me point out that there exist Calabi-Yau varieties $X$ which are not Frobenius split.  For such $X$, by Theorem~\ref{th:fsplit}, $B^1_X$ is not an ACM bundle. For example consider the Fermat quartic surface in $\P^3$ given by:
$$x^4+y^4+z^4+w^4=0.$$
For $p\congruent3\bmod{4}$ this is not Frobenius split (in fact $X$ is a supersingular $K3$ surface).
\erem

\begin{remark}\label{re:quartic}
Preceding results lead us to the following questions:
\benum[label={\bf (\arabic{*})}]
\item\label{qu:acm} What assumptions on $X\into \P^n$ are sufficient to ensure $B^1_X$ is an ACM bundle?
\item\label{qu:acm2} What conditions on $X$ are sufficient to ensure that $F_*(\O_X)$ is an ACM bundle?
\eenum
If $\dim(X)\geq 2$ then any Ulrich bundle is an ACM bundle and also satisfies a strong restriction on its Hilbert polynomial (see \cite{coskun13}). It appears that this Hilbert polynomial condition is the one which fails to hold for $B^1_X$ (and its twists) because this Hilbert polynomial condition is rather restrictive. For example this already is the case for smooth quartic surfaces in $\P^3$ (this is an easy calculation using (\cite{coskun13})).
On the other hand one sees from preceding results that there are interesting classes of smooth, projective variety for which \ref{qu:acm} and \ref{qu:acm2} have an affirmative answer. But I do not believe that the list of such varieties  which I provide here is exhaustive. 
\end{remark}

%
\bibliographystyle{plain}
\bibliography{ulrich,../../master/joshi,../../master/master6}
\end{document}